\documentclass[11pt]{amsart}
\usepackage{xcolor,mathtools}
\allowdisplaybreaks

\include{package}
\definecolor{grn}{rgb}{0,0.6,0}
\definecolor{mrn}{rgb}{0.3,0,0}
\definecolor{blue}{rgb}{0,0,0.7}
\definecolor{Mygray}{rgb}{0.75,0.75,0.75}
\definecolor{auburn}{rgb}{0.43, 0.21, 0.1}
\definecolor{britishracinggreen}{rgb}{0.0, 0.26, 0.15}
\definecolor{taupe}{rgb}{0.28, 0.24, 0.2}
\usepackage{amsmath,amssymb}
\newtheorem{theorem}{Theorem}

\newtheorem{proposition}{Proposition}[section]

\newtheorem{remark}{Remark}

\newtheorem*{ack}{Acknowledgements}

\allowdisplaybreaks
\usepackage{latexsym,amssymb}
\usepackage{amssymb}
\usepackage{graphicx}
\usepackage{ textcomp }
\usepackage{tikz}
\usetikzlibrary{shapes.geometric, arrows}
\tikzstyle{startstop} = [rectangle, rounded corners, minimum width=3cm, minimum height=1cm,text centered, draw=black, fill=white!30]
\tikzstyle{io} = [trapezium, trapezium left angle=70, trapezium right angle=110, minimum width=3cm, minimum height=1cm, text centered, draw=black, fill=white!30]
\tikzstyle{process} = [rectangle, minimum width=2cm, minimum height=1cm, text centered, draw=black, fill=white!30]
\tikzstyle{decision} = [rectangle, minimum width=1cm, minimum height=1cm, text centered, draw=black, fill=white!30]
\tikzstyle{arrow} = [thick,->,>=stealth]
\usepackage[left=2cm,right=2cm,top=2cm,bottom=2cm]{geometry}

\allowdisplaybreaks

\begin{document}
\baselineskip=14.5pt
\title[Elliptic curves with rank at least $3$]{Two infinite families of elliptic curves with Mordell-Weil rank at least $3$}

\author{Pankaj Patel, Debopam Chakraborty and Jaitra Chattopadhyay}
\address[Pankaj Patel and Debopam Chakraborty]{Department of Mathematics, BITS-Pilani, Hyderabad campus, Hyderabad, INDIA}
\address[Jaitra Chattopadhyay]{Department of Mathematics, Siksha Bhavana, Visva-Bharati, Santiniketan - 731235, West Bengal, India}

\email[Pankaj Patel]{p20200452@hyderabad.bits-pilani.ac.in}

\email[Debopam Chakraborty]{debopam@hyderabad.bits-pilani.ac.in}

\email[Jaitra Chattopadhyay]{jaitra.chattopadhyay@visva-bharati.ac.in; chat.jaitra@gmail.com}

\begin{abstract}
In this paper, we consider two infinite parametric families of elliptic curves defined over $\mathbb{Q}$ given by the equations $E_{a,b} : y^{2} = x^{3} - a^{2}x + b^{2}$ and $E^{\prime}_{a,b} : y^{2} = x^{3} - a^{2}x + b^{6}$, where $a,b \in \mathbb{N}$ satisfy certain mild conditions.  We prove that the torsion group of $E_{a,b}(\mathbb{Q})$ is trivial and the Mordell-Weil ranks of both $E_{a,b}(\mathbb{Q})$ and  $E^{\prime}_{a,b}(\mathbb{Q})$ are at least $3$ for infinitely many choices of $a$ and $b$ by using the N\'{e}ron-Tate height of a rational point and by exploiting the unit group of the ring of integers of $\mathbb{Q}(\sqrt{3})$. This is an extension of the results of Brown-Myres and Fujita-Nara where lower bounds of the ranks were provided under the assumption that $a = 1$ or $b = 1$. Also, our families of elliptic curves vastly generalize the curves recently investigated by Hatley and Stack. 
\end{abstract}

\renewcommand{\thefootnote}{}



\footnote{\emph{We confirm that all the data are included in the article. There is no conflict of interest.}}

\renewcommand{\thefootnote}{\arabic{footnote}}
\setcounter{footnote}{0}

\maketitle

\section{introduction}

Let $E$ be an elliptic curve defined over $\mathbb{Q}$. Then the set of all  rational points on $E$ forms an abelian group, denoted by $E(\mathbb{Q})$. The celebrated Mordell-Weil theorem states that $E(\mathbb{Q})$ is finitely generated, and therefore by the structure theorem of finitely generated abelian groups, we have $E(\mathbb{Q}) \simeq \mathbb{Z}^{r}\oplus E(\mathbb{Q})_{\rm{tors}}$, where $r$ is the Mordell-Weil rank of $E(\mathbb{Q})$ and $E(\mathbb{Q})_{\rm{tors}}$ is the torsion subgroup of it. The Mordell-Weil rank $r$ of $E(\mathbb{Q})$ is an extremely important object of study in number theory since it is related to various seemingly disparate problems in mathematics. A classic instance of this is the congruent number problem that asks for a characterization of all positive integers $n$ that can be realized as the area of a right-angled rational-sided triangle. This is equivalent to the assertion that the {\it congruent number elliptic curve} $y^{2} = x^{3} - n^{2}x$ has positive Mordell-Weil rank. 

\smallskip

It is a folklore conjecture that almost all elliptic curves defined over $\mathbb{Q}$ have Mordell-Weil ranks either $0$ or $1$. Therefore, it is useful to find elliptic curves whose ranks are at least $2$. In \cite{bro}, Brown and Myres considered the elliptic curve $E_{m} : y^{2} = x^{3} - x + m^{2}$ and proved that ${\rm{rank}}(E_{m}) \geq 3$ for infinitely many values of $m$. Fujita and Nara \cite{nara} generalized this and considered the family of elliptic curves given by $E_{m,n} : y^{2}  = x^{3} - m^{2}x + n^{2}$. They proved some lower bounds of ${\rm{rank}}(E_{m,n})$ for some particular values of $m$ and $n$. Recently, Juyal and Kumar \cite{juy} proved that ${\rm{rank}}(E_{m,p}) \geq 2$, where $p$ is an odd prime number and $m \geq 1$ is an integer with $m \not\equiv 0 \pmod {3}$ and $m \equiv 2 \pmod {32}$. Later, Chakraborty and Sharma \cite{cha} proved that ${\rm{rank}}(E_{m,pq}) \geq 2$, where $m \geq 1$ is an integer with $m \not\equiv 0 \pmod {3}$, $m \not\equiv 0 \pmod {4}$, $m \equiv 2 \pmod {64}$ and $p$ and $q$ are distinct odd prime numbers not dividing $m$. Very recently, Ghosh \cite{arka} further extended the results for $E_{m,pqr}$ where $p,q,r$ are three distinct prime numbers. Many authors (cf. \cite{ant}, \cite{eki}, \cite{arka}, \cite{integers}, \cite{rout}, \cite{self-rama}, \cite{tad1}, \cite{tad2}) have considered other interesting families of elliptic curves and investigated their torsion subgroups, Mordell-Weil ranks and Selmer groups. 

\smallskip

In this article, we consider two infinite parametric familes of elliptic curves and study their torsion subgroups and Mordell-Weil ranks. This is a further generalization of the curves considered by Fujita and Nara \cite{nara}. More precisely, we prove the following two theorems.
\begin{theorem}\label{1st}
Let $a$ and $b$ be positive integers with $\gcd(a,b) = 1$ and $b$ is an odd integer. If $3 \mid b$ and $4 \nmid a$, then the torsion group of $E_{a,b} : y^{2} = x^{3} - a^{2}x + b^{2}$ is trivial. Moreover, there exist infinitely many pairs of positive integers $(a,b)$ with $a^{2} = 3b^{2} + 1$ such that the Mordell-Well rank of $E_{a,b} : y^{2} = x^{3} - a^{2}x + b^{2}$ is at least $3$ and $E_{a,b}(\mathbb{Q})_{\rm{tors}}$ is trivial. 
\end{theorem}
\begin{remark}
In Theorem \ref{1st}, if $b$ is the product of two distinct prime numbers, then we get back the results obtained in \cite{cha}. Also, if $b$ is the product of three distinct prime numbers, then we recover the results obtained in \cite{arka}. Thus Theorem \ref{1st} vastly generalizes the results of \cite{cha} and \cite{arka}.
\end{remark} 
\begin{theorem}\label{3rd}
For integers $a, b$ with $0 < a < b$, the Mordell-Weil rank of the elliptic curve $E^{\prime}_{a,b} : y^{2} = x^{3} - a^{2}x + b^{6}$ is at least $3$. 
\end{theorem}
\begin{remark}
Hatley and Stack \cite{integers} considered the elliptic curve $y^{2} = x^{3} - x + m^{6}$, where $m$ is an odd positive integer divisible by $3$, and proved that its Mordell-Weil rank is at least $3$. Therefore, Theorem \ref{3rd} clearly generalizes their result. 
\end{remark}

\section{N\'{e}ron-Tate height and the Proof of Theorem \ref{1st}}

Let $E_{a,b} : y^{2} = x^{3} - a^{2}x + b^{2}$, where $a,b \in \mathbb{N}$ are such that $\gcd(a,b) = 1$ and $b$ is an odd integer. First, we prove that $E_{a,b}(\mathbb{Q})$ has no $2$-torsion point. For that, if possible let us assume that $P \in E_{a,b}(\mathbb{Q})$ is a $2$-torsion point. Then by Lutz-Nagell theorem, it follows that $P = (n,0)$ for some $n \in \mathbb{Z}$. Since $P = (n,0) \in E_{a,b}(\mathbb{Q})$, we have $n^{3} - a^{2}n + b^{2} = 0$, which implies that $n(n - a)(n + a) = -b^{2}$.

\smallskip

We claim that either $v_{p}(n) = 0 \mbox{ or } v_{p}(n) = v_{p}(b^{2})$ for all prime divisor $p$ of $b$. For, if $p$ is a prime divisor of $b$ with $v_{p}(n) < v_{p}(b^{2})$, then $p^{v_{p}(b^{2})} \mid n(n - a)(n + a)$ implies that $p \mid (n - a)$ or $p \mid (n + a)$. This forces $p \mid a$ which is a contradiction to the hypothesis that $\gcd(a,b) = 1$. Thus $v_{p}(n) = 0 \mbox{ or } v_{p}(b^{2})$. 

\smallskip

\noindent
{\bf Case I.} Assume that $n > 0$. Then $n + a > 0$ and from $n(n + a)(n - a) = -b^{2}$, it follows that $n - a < 0$. Hence $a - n > 0$. Since $\gcd(a,b) = 1$ and $n(n + a) \mid b^{2}$, we conclude that $a, a + n$ and $n - a$ are pairwise relatively prime. Consequently, we have $n + a = u^{2}$ and $a - n = v^{2}$ for some $u,v \in \mathbb{Z}$. Moreover, $b$ is odd implies that $n$ and $n + a$ are both odd integers and thus $a$ is even. From this, we obtain that both $u$ and $v$ are odd integers. Therefore, $$2a = (a + n) + (a - n) = u^{2} + v^{2} \equiv 1 + 1 \equiv 2 \pmod {4},$$ which is impossible because $a$ is even implies that $2a \equiv 0 \pmod {4}$. 

\smallskip

\noindent
{\bf Case II.} Assume that $n < 0$. Then $n(n + a)(n - a) = -b^{2}$ implies that $(n + a)(n - a) > 0$. Since $a > 0$, we must have $n + a < 0$ and $n - a < 0$. Again, the coprimality of $n + a$ and $n - a$ implies that $n + a = -r^{2}$ and $n - a = -s^{2}$ for some odd integers $r$ and $s$. Consequently, we obtain $$2a = (n + a) - (n - a) = s^{2} - r^{2} \equiv 0 \pmod {8},$$ which forces that $a \equiv 0 \pmod {4}$, a contradiction to our hypothesis. We therefore conclude that there is no $2$-torsion points in $E_{a,b}(\mathbb{Q})$. 

%
%
Now, we prove that $E_{a,b}(\mathbb{Q})$ has no element of finite order. We observe that the discriminant $\Delta(E_{a,b})$ of $E_{a,b}$ is a divisor of $-16(-4a^{6} + 27b^{4})$. Since $3 \mid b$, we have that $-16(-4a^{6} + 27b^{4}) \equiv 64a^{6} \pmod {3}$. Since $\gcd(a,b) = 1$ and $3$ divides $b$, we conclude that $\gcd(3,a) = 1$. Consequently, $\Delta(E_{a,b}) \not\equiv 0 \pmod {3}$ and thus $3$ is a prime of good reduction for $E_{a,b}$. Using $a^{2} \equiv 1 \pmod {3}$, we obtain $$\bar{E}_{a,b} : y^{2} \equiv x^{3} - x \pmod {3}.$$ It can be directly checked that $\bar{E}_{a,b}(\mathbb{F}_{3}) = \{\mathcal{O}, (\bar{0},\bar{0}), (\bar{1},\bar{0}), (\bar{2},\bar{0})\}$. That is, $\bar{E}_{a,b}(\mathbb{F}_{3})$ is a group of order $4$. Now, since $3$ is a prime of good reduction, there exists an injective group homomorphism from $E_{a,b}(\mathbb{Q})_{\rm{tors}}$ to $\bar{E}_{a,b}(\mathbb{F}_{3})$. It then follows that either $E_{a,b}(\mathbb{Q})_{\rm{tors}} = \{\mathcal{O}\}$ or it has an element of order $2$. Since we have already proved that $E_{a,b}(\mathbb{Q})$ has no $2$-torsion point, we conclude that $E_{a,b}(\mathbb{Q})_{\rm{tors}} = \{\mathcal{O}\}$. 

%
%
%


Next, we prove that the Mordell-Weil rank of the elliptic curve $E_{a,b} : y^{2} = x^{3} - a^{2}x + b^{2}$, where $a$ and $b$ are integers with $a^{2} = 3b^{2} + 1$, is at least $3$ for all positive integers $a$ and $b$. Then we prove that we can force the conditions $a$ is even, $b$ is odd with $4 \nmid a$ and $3 \mid b$ for infinitely many of them so that the torsion group turns out to be trivial. We note that $P_{1} = (0,b)$, $P_{2} = (a,b)$ and $P_{3} = (-1,2b)$ are three rational points on $E_{a,b}$. We shall prove that these three points are $\mathbb{Z}$-linearly independent. We accomplish this by using the canonical height of $P_{1}$, $P_{2}$ and $P_{3}$ and establishing that the height matrix corresponding to this is non-singular. We closely follow the treatment that is provided in \cite{integers}.  

\smallskip

We recall that for a rational number $\alpha = \frac{p}{q}$, the height of $\alpha$ is defined by $h(\alpha) = \log(\max\{|p|,|q|\})$. Using this a height function $H$ on the set of all rational points of $E_{a,b}$ is defined by setting $H(P) = h(x(P))$, where $P = (x(P),y(P)) \in E_{a,b}(\mathbb{Q})$. Finally, we define the N\'{e}ron-Tate height (also known as the  canonical height) $\hat{h}$ of $P$ by letting $\hat{h}(P) = \frac{1}{2} \displaystyle\lim_{N \to \infty}\frac{H(2^{N}P)}{4^{N}}$, because it is a well-known result in the arithmetic of elliptic curves that the limit always exists. Below, we list a few results that are useful to prove that ${\rm{rank}}(E_{a,b}) \geq 3$. 

\begin{proposition} \cite[Proposition 8]{integers}\label{integer-proposition}
Let $E$ be an elliptic curve defined over $\mathbb{Q}$ and let $\hat{h}$ be the N\'{e}ron-Tate height function. Then the following hold true.
\begin{enumerate}
\item We have $\hat{h}(P) \geq 0$ for every $P \in E(\mathbb{Q})$. Also, equality holds if and only if $P \in E(\mathbb{Q})_{{\rm{tors}}}$. 

\medskip

\item We have $\hat{h}(P_{1} + P_{2}) + \hat{h}(P_{1} - P_{2}) = 2\hat{h}(P_{1}) + 2\hat{h}(P_{2})$ for every $P_{1}, P_{2} \in E(\mathbb{Q})$. This is sometimes referred to as the parallelogram law.

\medskip

\item We have $\hat{h}(tP) = t^{2}\hat{h}(P)$ for every $P \in E(\mathbb{Q})$ and $t \in \mathbb{Z}$,.
\end{enumerate}
\end{proposition}

The N\'{e}ron-Tate height pairing is defined as the map from $E(\mathbb{Q}) \times E(\mathbb{Q}) \to \mathbb{R}$ by the equation $\langle P_{1},P_{2} \rangle = \hat{h}(P_{1} + P_{2}) - \hat{h}(P_{1}) - \hat{h}(P_{2})$. It is a symmetric bilinear map and $P_{1},\ldots,P_{k} \in E(\mathbb{Q})$ are $\mathbb{Z}$-linearly independent if and only if the height-pairing $k \times k$-matrix $[\langle P_{i},P_{j}\rangle]_{1 \leq i,j \leq k}$ is non-singular. 

\smallskip

From the definition of the N\'{e}ron-Tate height, it is evident that we need to calculate the $x$-co-ordinate of the points $2^{N}P$. If $P = (x,y) \in E_{a,b}(\mathbb{Q})$, then using the duplication formula we obtain $$x(2P) = \dfrac{x^{4} + 2a^{2}x^{2} - 8b^{2}x + a^{4}}{4(x^{3} - a^{2}x + b^{2})}.$$ Using the relation $a^{2} = 3b^{2} + 1$, we can rewrite this as $$x(2P) = \dfrac{x^{4} + 2a^{2}x^{2} + a^{4} - \frac{8}{3}a^{2}x + \frac{8}{3}x}{4x^{3} - 4a^{2}x + \frac{4}{3}a^{2} - \frac{4}{3}}.$$ We observe that $x(2P_{1}) = \frac{(3b^{2} + 1)^{2}}{4b^{2}} = \frac{a^{4}}{\frac{4}{3}a^{2} - \frac{4}{3}}$ and the exponent of $a$ in the sequence $\{x(2^{N}P)\}_{N \geq 1}$ grows as $4, 16, 64, 256,\ldots$. More precisely, using induction on $N$, we prove that $x(2^{N}P_{1}) = \dfrac{a^{4^{N}} + f(a)}{g(a)}$, where $f(a)$ and $g(a)$ are polynomial expressions in $a$ with rational coefficients of degrees $4^{N} - 2$ in $a$ each for any integer $N \geq 2$. 


\smallskip

We see by the duplication formula that $x(2P_{1}) = \dfrac{a^{4}}{\frac{4}{3}a^{2} - \frac{4}{3}}$ and hence $x(4P_{1}) = \dfrac{u_{2,a}}{v_{2,a}}$, where
\begin{eqnarray*}
u_{2,a} &=& a^{16} + 2a^{10}\left(\frac{4}{3}a^{2} - \frac{4}{3}\right)^{2} + a^{4}\left(\frac{4}{3}a^{2} - \frac{4}{3}\right)^{4} -\frac{8}{3}a^{6}\left(\frac{4}{3}a^{2} - \frac{4}{3}\right)^{3} \\ &+& \frac{8}{3}a^{4}\left(\frac{4}{3}a^{2} - \frac{4}{3}\right)^{3}
\end{eqnarray*}

and

\begin{eqnarray*}
v_{2,a} &=& 4a^{12}\left(\frac{4}{3}a^{2} - \frac{4}{3}\right) - 4a^{6}\left(\frac{4}{3}a^{2} - \frac{4}{3}\right)^{3} + \frac{4}{3}a^{2}\left(\frac{4}{3}a^{2} - \frac{4}{3}\right)^{4} \\ &-& \frac{4}{3}\left(\frac{4}{3}a^{2} - \frac{4}{3}\right)^{4}.
\end{eqnarray*}

We immediately see that the statement holds true for $N = 2$ since $u_{2,a}$ is a polynomial expression in $a$ with rational coefficients of degree $16$. Now, we assume that $x(2^{N}P_{1}) = \dfrac{a^{4^{N}} + f(a)}{g(a)}$, where $f(a)$ and $g(a)$ are polynomial expressions in $a$ with rational coefficients of degrees $4^{N} - 2$ in $a$ each for some integer $N \geq 2$. Then again using the duplication formula, we obtain $x(2^{N + 1}P_{1}) = x(2\cdot 2^{N}P_{1}) = \dfrac{u_{a}}{v_{a}}$, where 

\begin{eqnarray*}
u_{a} &=& (a^{4^{N}} + f(a))^{4} + 2a^{2}(a^{4^{N}} + f(a))^{2}g(a)^{2} - \frac{8}{3}a^{2}(a^{4^{N}} + f(a))g(a)^{3} \\ &+& \frac{8}{3}(a^{4^{N}} + f(a))g(a)^{3} + a^{4}g(a)^{4} 
\end{eqnarray*}

and 

\begin{eqnarray*}
v_{a} &=& 4(a^{4^{N}} + f(a))^{3}g(a) - 4a^{2}(a^{4^{N}} + f(a))g(a)^{3} + \frac{4}{3}a^{2}g(a)^{4} - \frac{4}{3}g(a)^{4}.
\end{eqnarray*}
We see that the first summand of $u_{a}$ has degree $4^{N + 1}$ in $a$ and using the facts that $f(a)$ and $g(a)$ have degrees $4^{N} - 2$ each in $a$, we conclude that the remaining summands of $u_{a}$ have respective degrees $4^{N + 1} - 2$, $4^{N + 1} - 4$, $4^{N + 1} - 6$ and $4^{N + 1} - 4$. Similarly, considering the successive summands of $v_{a}$, we see that the respective degrees in $a$ are $4^{N + 1} - 2$, $4^{N + 1} - 4$, $4^{N + 1} - 6$ and $4^{N + 1} - 8$. Therefore, by induction of $N$, we conclude that $H(2^{N}P_{1}) = \log(a^{4^{N}} + g_{N}(a))$, where $g_{N}(a)$ is a polynomial in $a$ of degree $4^{N} - 2$. Therefore, the N\'{e}ron-Tate height of $P_{1}$ is $$\hat{h}(P_{1}) = \frac{1}{2}\displaystyle\lim_{N \to \infty}\dfrac{H(2^{N}P_{1})}{4^{N}} = \frac{1}{2}\displaystyle\lim_{N \to \infty}\dfrac{\log(a^{4^{N}} + g_{N}(a))}{4^{N}} = \dfrac{1}{2}.$$
 
It is worth mentioning that without loss of any generality, we have been considering the logarithm to the base $a$.

\smallskip

Similar calculations for $P_{2}$ and $P_{3}$ show that $\hat{h}(P_{2}) = \hat{h}(P_{3}) = \hat{h}(P_{1} + P_{2}) = \frac{1}{2}$. Using the parallelogram law from Proposition \ref{integer-proposition}, we obtain $$\langle P_{1}, P_{2} \rangle = \hat{h}(P_{1} + P_{2}) - \hat{h}(P_{1}) - \hat{h}(P_{2}) = \frac{1}{2} - \frac{1}{2} - \frac{1}{2} = -\frac{1}{2}.$$ The parallelogram law applied to $P_{1}, P_{3}$ and $P_{2},P_{3}$ yield $\langle P_{1},P_{3} \rangle = 0 =  \langle P_{2},P_{3} \rangle$. Consequently, the height-pairing matrix corresponding to $P_{1}, P_{2}$ and $P_{3}$ is $\begin{bmatrix}
1 & -\frac{1}{2} & 0 \\
-\frac{1}{2} & 1 & 0 \\
0 & 0 & 1
\end{bmatrix}$. Since this is a non-singular matrix, we conclude that $P_{1}, P_{2}$ and $P_{3}$ are $\mathbb{Z}$-linearly independent points on $E_{a,b}(\mathbb{Q})$. Hence ${\rm{rank}}(E_{a,b}(\mathbb{Q})) \geq 3$. 

\smallskip

Now, it only remains to prove that there exist infinitely many positive integers $a$ and $b$ satisfying $a^{2} - 3b^{2} = 1$, $b$ is odd, $3 \mid b$ and $4 \nmid a$. We observe that if $a,b \geq 1$ are integers such that $a^{2} - 3b^{2} = 1$, then $a + b\sqrt{3}$ is a unit in the ring $\mathbb{Z}[\sqrt{3}]$, which is the ring of integers of the real quadratic field $\mathbb{Q}(\sqrt{3})$. Also, from the equation $a^{2} - 3b^{2} = 1$, it directly follows that $\gcd(a,b) = 1$. 

\smallskip

We note that $\varepsilon = 2 + \sqrt{3}$ is a unit of $\mathbb{Z}[\sqrt{3}]$ and therefore, $\varepsilon^{n}$ is also a unit for any integer $n \geq 1$. Let $\varepsilon^{n} = (2 + \sqrt{3})^{n} = a_{n} + b_{n}\sqrt{3}$, where $a_{n}, b_{n} \in \mathbb{Z}$. Then from the equation $a_{n + 1} + b_{n + 1}\sqrt{3} = (2 + \sqrt{3})(a_{n} + b_{n}\sqrt{3})$, we obtain the recurrence relations
\begin{equation}\label{recurrence}
a_{n + 1} = 2a_{n} + 3b_{n} \qquad \mbox{ and } \qquad b_{n + 1} = a_{n} + 2b_{n},
\end{equation}
with the initial conditions $a_{1} = 2$ and $b_{1} = 1$. Now, $a_{n}^{2} - 3b_{n}^{2} = 1$ implies that $a_{n}$ and $b_{n}$ are of different parity. Therefore, if $a_{n}$ is even, then $b_{n}$ is odd and consequently, we have $b_{n}^{2} \equiv 1 \pmod {8}$. Therefore, from $a_{n}^{2} - 3b_{n}^{2} = 1$, it follows that $a_{n}^{2} \equiv 4 \pmod {8}$ and hence $4 \nmid a_{n}$.

\smallskip

Next, we establish that if $3 \mid b_{n}$, then $3 \mid b_{n + 3}$. Clearly, if $3 \mid b_{n}$, then $3 \nmid a_{n}$ because $a_{n}^{2} - 3b_{n}^{2} = 1$. Therefore, either $a_{n} \equiv 1 \pmod {3}$ or $a_{n} \equiv 2 \pmod {3}$. We first consider the case where $a_{n} \equiv 1 \pmod {3}$. In that case, we have $a_{n + 1} = 2a_{n} + 3b_{n} \equiv 2 \pmod {3}$ and $b_{n + 1} = a_{n} + 2b_{n} \equiv 1 \pmod {3}$. It then follows that $a_{n + 2} = 2a_{n + 1} + 3b_{n + 1} \equiv 2\cdot 2 \equiv 1 \pmod {3}$ and $b_{n + 2} = a_{n + 1} + 2b_{n + 1} \equiv 2 + 2\cdot 1 \equiv 1 \pmod {3}$. Consequently, we have $b_{n + 3} = a_{n + 2} + 2b_{n + 2} \equiv 1 + 2\cdot 1 \equiv 0 \pmod {3}$. 

\smallskip

Now, we consider the remaining case where $a_{n} \equiv 2 \pmod {3}$. In this case, we have $a_{n + 1} = 2a_{n} + 3b_{n} \equiv 1 \pmod {3}$ and $b_{n + 1} = a_{n} + 2b_{n} \equiv 2 \pmod {3}$. It then follows that $a_{n + 2} = 2a_{n + 1} + 3b_{n + 1} \equiv 2 \pmod {3}$ and $b_{n + 2} = a_{n + 1} + 2b_{n + 1} \equiv 2 \pmod {3}$. Hence $b_{n + 3} = a_{n + 2} + 2b_{n + 2} \equiv 2 + 2\cdot 2 \equiv 0 \pmod {3}$. Thus we established that if $3 \mid  b_{n}$, then $3 \mid b_{n + 3\ell}$ for all integer $\ell \geq 0$. 

\smallskip

From the recurrence relations in \eqref{recurrence}, we see that if $b_{n}$ is odd, then so is $b_{n + 2k}$ for every integer $k \geq 0$. Now, $\varepsilon^{3} = 26 + 15\sqrt{3}$ shows that $b_{3} = 15$ is odd and is divisible by $3$. Hence $b_{3 + 6\ell}$ is odd and divisible by $3$ for all integer $\ell \geq 0$. Thus we obtain infinitely many positive integers $a$ and $b$ satisfying $a^{2} - 3b^{2} = 1$, $b$ is odd, $3 \mid b$ and $4 \nmid a$. This completes the proof of the theorem. $\hfill\Box$

\section{Proof of Theorem \ref{3rd}}

In this section, we deal with the elliptic curve $E^{\prime}_{a,b} : y^{2} = x^{3} - a^{2}x + b^{6}$ and prove that its Mordell-Weil rank is at least $3$. Like the proof of Theorem \ref{1st}, here also we consider the N\'{e}ron-Tate height of three particular rational points on $E^{\prime}_{a,b}(\mathbb{Q})$ and prove that the height-pairing matrix is non-singular. This will prove that the Mordell-Weil rank of $E^{\prime}_{a,b(\mathbb{Q})}$ is at least $3$.

\smallskip

We can directly check form the defining equation of $E^{\prime}_{a,b}$ that $P_{1} = (0,b^{3})$, $P_{2} = (a,b^{3})$ and $P_{3} = (-b^{2},ab)$ are elements of $E^{\prime}_{a,b}(\mathbb{Q})$. We shall prove that $P_{1}$, $P_{2}$ and $P_{3}$ are $\mathbb{Z}$-linearly independent and thus it will be proved that ${\rm{rank}}(E^{\prime}_{a,b}(\mathbb{Q})) \geq 3$. We follow exactly the method adopted in the proof of Theorem \ref{1st} and compute the N\'{e}ron-Tate height of these points. Then the non-singularity of the height-pairing matrix yields that the points $P_{1}$, $P_{2}$ and $P_{3}$ are $\mathbb{Z}$-linearly independent. 

\smallskip

Note that if $P = (x,y)$ is an element of $E^{\prime}_{a,b}(\mathbb{Q})$, then the duplication formula gives $$x(2P) = \dfrac{x^{4} + 2a^{2}x^{2} - 8b^{6}x + a^{4}}{4x^{3} - 4a^{2}x + 4b^{6}}.$$ Therefore, for the point $P_{1} = (0,b^{3})$, we have $x(2P_{1}) = \frac{a^{4}}{4b^{6}}$. Again, we have
\begin{equation}\label{th-2-x-4p}
x(4P_{1}) = \dfrac{a^{16} + 32a^{10}b^{12} - 256a^{4}b^{24}}{1024b^{30} - 256a^{6}b^{18} + 16a^{12}b^{6}}.
\end{equation}
We claim that $x(2^{N}P_{1}) = \dfrac{f(a,b) + f_{N}(a,b)}{g(a,b) + g_{N}(a,b)}$ for every integer $N \geq 2$, where $f(a,b), g(a,b), f_{N}(a,b), g_{N}(a,b)$ are polynomial expressions in $a$ and $b$, the exponents of $b$ in $f(a,b)$ and $g(a,b)$ are $2\cdot 4^{N} - 8$ and $2\cdot 4^{N} - 2$, respectively. Moreover, the exponent of $b$ in $f_{N}(a,b)$ and $g_{N}(a,b)$ are strictly smaller than those in $f(a,b)$ and $g(a,b)$, respectively. 

\smallskip

We immdeiately see from the expression of $x(4P_{1})$ in equation \eqref{th-2-x-4p} that the claim holds true for $N = 2$. Now, assume that the claim holds true for some integer $N \geq 2$. Then we have $x(2^{N + 1}P_{1}) = x(2\cdot 2^{N}P_{1}) = \dfrac{u(a,b)}{v(a,b)}$, where 
\begin{eqnarray*}
u(a,b) &=& (f(a,b) + f_{N}(a,b))^{4} + 2a^{2}(f(a,b) + f_{N}(a,b))^{2}(g(a,b) + g_{N}(a,b))^{2} \\ &-& 8b^{6}(f(a,b) + f_{N}(a,b))(g(a,b) + g_{N}(a,b))^{3} + a^{4}(g(a,b) + g_{N}(a,b))^{4}
\end{eqnarray*}

and 

\begin{eqnarray*}
v(a,b) &=& 4(f(a,b) + f_{N}(a,b))^{3}(g(a,b) + g_{N}(a,b)) \\ &-& 4a^{2}(f(a,b) + f_{N}(a,b))(g(a,b) + g_{N}(a,b))^{3} + 4b^{6}(g(a,b) + g_{N}(a,b))^{4}. 
\end{eqnarray*}

Hence we obtain that the exponents of $b$ in the successive summands of $u(a,b)$ are $2\cdot 4^{N + 1} - 32$, $2\cdot 4^{N + 1} - 20$, $2\cdot 4^{N + 1} - 8$ and $2\cdot 4^{N + 1} - 8$, respectively. Also, the exponents of $b$ in the successive summands of $v(a,b)$ are $2\cdot 4^{N + 1} - 26$, $2\cdot 4^{N + 1} - 14$ and $2\cdot 4^{N + 1} - 2$, respectively. Consequently, the claim is established by induction on $N$. Hence by considering logarithms to the base $b$, for any integer $N \geq 2$, we have $$\hat{h}(P_{1}) = \dfrac{1}{2}\displaystyle\lim_{N \to \infty}\dfrac{H(2^{N}P_{1})}{4^{N}} = \dfrac{1}{2}\displaystyle\lim_{N \to \infty}\dfrac{2\cdot 4^{N} - 2}{4^{N}} = 1.$$ 

Hence we have $\langle P_{1},P_{1} \rangle = \hat{h}(P_1 + P_{1}) - \hat{h}(P_{1}) - \hat{h}(P_{1}) = \hat{h}(2P_{1}) - 2\hat{h}(P_{1}) = 4\hat{h}(P_{1}) - 2 = 2$. Also, working with the points $P_{2}$ and $P_{3}$, we obtain $\hat{h}(P_{2}) = \hat{h}(P_{3}) = 1$. Moreover, using Proposition \ref{integer-proposition}, we obtain $$\langle P_{1},P_{2} \rangle = -1, \qquad \langle P_{1},P_{3} \rangle = 0 \qquad \mbox{ and } \qquad \langle P_{2},P_{3} \rangle = 1.$$ Consequently, the N\'{e}ron-Tate height-pairing matrix is $\begin{bmatrix}
2 & -1 & 0 \\
-1 & 2 & 1 \\
0 & 1 & 2
\end{bmatrix}$ which is non-singular. Hence ${\rm{rank}}(E^{\prime}_{a,b}(\mathbb{Q})) \geq 3$. This completes the proof of Theorem \ref{3rd}. $\hfill\Box$ 

\begin{ack}
We are grateful to our respective institutes for providing excellent facilities to carry out this research. The first author is supported by UGC (Grant no. 211610060298) and he thankfully acknowledges the same.
\end{ack}

\end{document}